\documentclass[12pt]{amsart}
\advance\textheight\topskip
\usepackage{times}
\allowdisplaybreaks

\usepackage[matrix,arrow,curve,dvips]{xy} 
\usepackage[dvips, bookmarks=false, hyperfigures=false,
setpagesize=false]{hyperref}

\usepackage[mathscr]{eucal}
\usepackage{mathptmx}

\newcommand{\tensor}{\mathbin{\otimes}}

\newcommand{\WW}[1]{\mathscr{W}^{(2)}_{#1}}
\newcommand{\ord}{{\operatorname{ord}}}

\newcommand{\Dynkin}[3]{\xymatrix@C28pt@1{\circ\ar@{-}[r]^(.55){#2}\ar@{}^{#1}[]&\circ\ar@{}^{#3}[]}}

\newcommand{\mfrac}[2]{\raisebox{.8pt}{\mbox{\small$\displaystyle\frac{#1}{#2}$}}}
\newcommand{\ffrac}[2]{\raisebox{.5pt}{\mbox{\footnotesize$\displaystyle\frac{#1}{#2}$}}}
\newcommand{\fffrac}[2]{\raisebox{.9pt}{\mbox{\scriptsize$\displaystyle\frac{#1}{#2}$}}}
\newcommand{\half}{%
  \mathchoice{\ffrac{1}{2}}{\frac{1}{2}}{\frac{1}{2}}{\frac{1}{2}}}

\newcommand{\hSL}[1]{\widehat{s\ell}(#1)}
\newcommand{\heis}{\widehat{\mathfrak{h}}}

\newcommand{\W}{\mathscr{W}}
\newcommand{\Wa}{\mathscr{W}_{\alpha}}
\newcommand{\Wb}{\mathscr{W}_{\beta}}
\newcommand{\Wba}{\mathscr{W}_{\beta\alpha}}
\newcommand{\Wab}{\mathscr{W}_{\alpha\beta}}
\newcommand{\Waba}{\mathscr{W}_{\alpha\beta\alpha}}
\newcommand{\Wbab}{\mathscr{W}_{\beta\alpha\beta}}
\newcommand{\Waabb}{\mathscr{W}_{\alpha\alpha\beta\beta}}

\newcommand{\Ea}{\mathscr{E}_{\alpha}}
\newcommand{\Eb}{\mathscr{E}_{\beta}}

\newcommand{\jplus}{j^+}
\newcommand{\jminus}{j^-}

\newcommand{\Jnaught}{J^0}

\newcommand{\oC}{\mathbb{C}}
\newcommand{\oZ}{\mathbb{Z}}
\newcommand{\Hel}[2]{{\small\textbf{#1}(#2)$_{\text{\cite{[Hel]}}}$}}

\newcommand{\G}[2]{\alpha_{#1}.\alpha_{#2}}
\newcommand{\Nich}{\mathfrak{B}}
\newcommand{\BX}{\Nich(X)}

\begin{document}

\title{Virasoro central charges for Nichols algebras}

\author[Semikhatov]{AM Semikhatov}

\address{Lebedev Physics Institute, Moscow 119991, Russia}

\begin{abstract}
  A Virasoro central charge can be associated with each Nichols
  algebra with diagonal braiding in a way that is invariant under the
  Weyl groupoid action.  The central charge takes very suggestive
  values for some items in Heckenberger's list of rank-$2$ Nichols
  algebras.  In particular, this might be viewed as an indication of
  the existence of reasonable logarithmic extensions of $W_3\equiv
  W\!A_2$, $W\!B_2$, and $W\!G_2$ models of conformal field theory.
  In the $W_3$ case, the construction of an \textit{octuplet} extended
  algebra---a counterpart of the triplet $(1,p)$ algebra---is
  outlined.
\end{abstract}

\maketitle%
\addtocounter{section}{-1}%
\thispagestyle{empty}%

\section{Introduction}
In~\cite{[STbr]}, we described a paradigm treating screening operators
in two-dimensional conformal field theory as a braided Hopf algebra, a
Nichols
algebra~\cite{[Nich],[AG],[AS-pointed],[Andr-remarks],[AS-onthe],
  [Heck-class],[AHS],[ARS]}.  This immediately suggests that the
inverse relation may also exist.  Is any finite-dimensional Nichols
algebra with diagonal braiding an algebra of screenings in some
conformal model?  This is a fascinating problem, especially
considering the recent remarkable development in the theory of Nichols
algebras---originally a ``technicality'' in Andruskiewitsch and
Schneider's program of classification of pointed Hopf algebras, which
has grown into a beautiful theory in and of itself (in addition to the
papers cited above and the references therein, also
see~\cite{[Heck-Weyl],[GHV],[GH-lyndon],[AHS],[AFGV],[AAY],
  [Ag-0804-standard],[Ag-1008-presentation],[Ag-1104-diagonal]}).
Diagonal braiding is assumed in what follows.
 
As many ``inverse'' problems, that of identifying a conformal field
model behind a given Nichols algebra is not necessarily well defined.
It is of course well known that screenings can be used to define
models of conformal field theory; in particular, defining logarithmic
models as kernels of screening~\cite{[FHST]} turned out to be
especially useful.  But passing from a Nichols algebras to screenings
involves various ambiguities.  Nevertheless, the central charges
associated with Nichols algebras in what follows have the nice
property of being invariant under the Weyl groupoid action---the
natural ``symmetry'' up to which Nichols algebras are
classified~\cite{[Heck-Weyl],[HY],[HS]}.

To proceed beyond the central charge identification, I restrict myself
to Nichols algebras of rank two (already a fairly large number in
terms of the possible conformal models).  All of these were listed by
Heckenberger~\cite{[Heck-1+1]} (the general classification, for any
rank, was achieved in~\cite{[Heck-class]} and was reproduced in a
different and independent way
in~\cite{[Ag-0804-standard],[Ag-1008-presentation],[Ag-1104-diagonal]}).
These notes\pagebreak[3] are in fact a compilation of the original
Heckenberger's list with explicit results on the presentation of some
Nichols algebras (obtained in~\cite{[Ag-0804-standard]} for the
standard type and in~\cite{[Hel]} in several nonstandard cases), and
with several CFT{} constructions added.  The extended algebra of a
logarithmic model---the octuplet algebra extending the $W_3$
algebra---is offered in only one case; the other CFT{} constructions
are merely a starting point for finding extended algebras.

For the Nichols algebra $\Nich(X)$ of a $\theta$-dimensional braided
linear space $(X,\Psi)$ (a rank-$\theta$ Nichols algebra), we fix a
basis $(F_i)_{1\leq i\leq\theta}$ in $X$ that diagonalizes the
braiding $\Psi:X\tensor X\to X\tensor X$,
\begin{equation}\label{q-matrix}
  \Psi(F_i\tensor F_j)=q_{i,j} F_j\tensor F_i,
\end{equation}
and call the matrix $(q_{i,j})_{1\leq i,j\leq\theta}$ the braiding
matrix.

The relation to conformal field theory is based on representing the
$F_i$ as screening operators
\begin{equation}\label{FF}
  F_i = \oint e^{\alpha_i.\varphi}
\end{equation}
acting in a space of bosonic fields (or simply ``bosons'').  Here,
$\varphi(z)=(\varphi^1(z),\dots,\varphi^{\theta}(z))$ is a
$\theta$-component boson field with OPEs~\eqref{ff-OPE}
(Appendix~\ref{app:vir}), the dot denotes Euclidean scalar product,
and $\alpha_i\in\oC^{\theta}$ are such that the screenings have the
self-braidings and the monodromies coincident with those
in~\eqref{q-matrix}:\footnote{Notably, the conditions on the braiding
  matrix elements selecting a Nichols algebra involve only the
  self-braidings $q_{i,i}$ and the \textit{monodromies}
  $q_{i,j}q_{j,i}$ for $i\neq j$.}
\begin{equation}\label{correspondence}
  \begin{aligned}
    e^{i\pi \alpha_i.\alpha_i}&=q_{i,i},\\
    e^{2i\pi \alpha_i.\alpha_j}&=q_{i,j}q_{j,i},\quad i\neq j.
  \end{aligned}
\end{equation}

The ambiguities inherent in passing from a braiding matrix to
screenings realized in terms of free bosons are numerous.  Already the
``$\theta$-boson space'' on which the $F_{i}$ act can be chosen
differently, e.g., including or not including exponentials
$e^{\omega.\varphi(z)}$, where $\omega$ ranges a lattice
in~$\oC^{\theta}$.  Furthermore, solving
relations~\eqref{correspondence} for $\alpha_i\in\oC^{\theta}$
involves taking logarithms, which introduces arbitrary integer
parameters.

And yet the idea to look for a conformal model corresponding to a
given Nichols algebra is not altogether meaningless
because \textit{the Virasoro central charge is invariant under the
  Weyl groupoid action}.  I go into some detail here because the
statement implicitly refers to a \textit{procedure} to deal with the
ambiguities such that the invariance is nevertheless ensured.

Whenever the $\alpha_i$ are linearly independent, the screenings
in~\eqref{FF} uniquely define a Virasoro algebra in their centralizer
in the space of differential polynomials in the $\partial\varphi^j(z)$
\ ($\partial=\partial/\partial z$). \ The general case is considered
in Appendix~\ref{app:vir}, and for $\theta=2$, for example, the
central charge of Virasoro algebra is
\begin{equation}\label{c-rank2}
  c=2-3\mfrac{\bigl(4 +  (\alpha_1.\alpha_1) (\alpha_2.\alpha_2)\bigr)
    (\alpha_1 - \alpha_2).(\alpha_1 - \alpha_2)
    + 4 (\alpha_1 - \alpha_2).\bigl((\alpha_1.\alpha_1) \alpha_2
    - (\alpha_2.\alpha_2) \alpha_1\bigr)}{(\alpha_1.\alpha_1)
    (\alpha_2.\alpha_2) - (\alpha_1.\alpha_2)^2}.
\end{equation}
On the Nichols algebra side, the Weyl groupoid action is defined as
follows~\cite{[Heck-Weyl],[HY],[HS],[CH]}.  There exists a generalized
Cartan matrix $(a_{i,j})_{1\leq i,j\leq \theta}$ such that $a_{i,i}=2$
and
\begin{equation}\label{cartan-test}
  q_{i,i}^{a_{i,j}} = q_{i,j}q_{j,i}\quad\text{or}\quad
  q_{i,i}^{1-a_{i,j}} = 1
\end{equation}
holds for each pair $i\neq j$.  The Weyl groupoid is generated by
pseudoreflections acting on the set of braiding matrices and defined
for any $k$, $1\leq k\leq\theta$.  The reflected braiding matrix has
the entries
\begin{equation}\label{refl-def}
  \mathfrak{R}^{(k)}(q_{i,j}) =
  q_{i,j} q_{i,k}^{-a_{k,j}} q_{k,j}^{-a_{k,i}} q_{k,k}^{a_{k,i}a_{k,j}}.
\end{equation}
It may or may not have the same generalized Cartan matrix.\footnote{If
  the diagonal braiding is of Cartan type, then Weyl reflections
  preserve the Cartan matrix.  If a generalized Cartan matrix (not of
  Cartan type) is the same for the entire class of Weyl-reflected
  braided matrices, then such a generalized Cartan matrix and the
  braiding matrix are said to belong to the \textit{standard} type.
  Nonstandard braidings do exist~\cite{[AA-gen],[Ag-0804-standard]}.}
The use of this tool has remarkably resulted in the classification of
Nichols algebras with diagonal braiding~\cite{[Heck-class]}.

With the screening momenta $\alpha_i\in\oC^{\theta}$, $1\leq
i\leq\theta$, defined such that~\eqref{correspondence} holds,
condition~\eqref{cartan-test} is ``lifted'' to the scalar products as
the condition
\begin{equation}\label{cartan-test-log}
  2\alpha_i.\alpha_j = a_{i,j} \alpha_i.\alpha_i
  \quad\text{or}\quad (1-a_{i,j})\alpha_i.\alpha_i = 2
\end{equation}
to be satisfied for each pair $i\neq j$.  Several particular choices
have been made in writing this (for example, the $2$, not some other
even integer, in the second relation).  

The Weyl reflections are now lifted to the scalar products by
``naively taking the logarithm'' of~\eqref{refl-def} (which amounts to
actual pseudoreflections in $\oC^\theta$):
\begin{equation}\label{Rk-act}
  \mathfrak{R}^{(k)}(\alpha_i.\alpha_j) =
  \alpha_i.\alpha_j - a_{k,j}\alpha_i.\alpha_k
  - a_{k,i}\alpha_k.\alpha_j
  + a_{k,i}a_{k,j} \alpha_k.\alpha_k.
\end{equation}
Weyl-reflecting the central charge amounts to replacing each
$\alpha_i.\alpha_j$ with $\mathfrak{R}^{(k)}(\alpha_i.\alpha_j)$ in
the system of equations that defines the central charge (see
Appendix~\ref{app:vir}).  

\smallskip

\noindent\textbf{Theorem.} \ \
\textit{The central charge of the Virasoro algebra centralizing
  $\theta$ screenings in the $\theta$-boson space is invariant
  under~\eqref{Rk-act} if conditions~\eqref{cartan-test-log} hold}.

\smallskip

This is proved in Appendix~\ref{app:vir}; for example, in
the rank-2 case (see~\eqref{c-rank2}), the Weyl-reflected central
charge is expressed rather explicitly as
\begin{multline*}
  \mathfrak{R}^{(1)}(c) - c = 
  \mfrac{3}{(\G{1}{1})(\G{2}{2}) - (\G{1}{2})^2}
  \bigl(2 \G{1}{2} - a_{1,2}\,\G{1}{1}\bigr)
  \bigl((a_{1,2}-1)\,\G{1}{1} + 2\bigr)
  \\
  \bigl(
  a_{1,2}^2\,\G{1}{1} 
  - a_{1,2}\,\G{1}{1} 
  - 2a_{1,2}\, \G{1}{2} 
  + 2 \G{2}{2}
  + 2 a_{1,2}
  - 4\bigr).
\end{multline*}
The product $\bigl(2 \alpha.\beta - a_{1,2}\,\alpha.\alpha\bigr)
\bigl((a_{1,2}-1)\,\alpha.\alpha + 2\bigr)$ vanishes
whenever~\eqref{cartan-test-log} holds for $i=1$ and $j=2$, and hence
$c$ is indeed invariant under $\mathfrak{R}^{(1)}$; the invariance
under $\mathfrak{R}^{(2)}$ is verified
similarly.

\subsection*{0.5. \ From Virasoro to extended algebras}
The central charge value alone does not specify a conformal field
theory uniquely.  In ``good'' cases, however---when the central charge
found from~\eqref{c-rank2} is a function of a (discrete)
parameter---the form of this dependence does suggest what type of
operators extend the Virasoro algebra and therefore what the resulting
conformal model is; and the centralizer of the screenings then turns
out to be sufficiently ample for an interesting conformal field theory
to live in it.  An exemplary case is the $W_3$ algebra, which
centralizes two screenings associated with
a braiding matrix such that $q_{1,1}=q_{2,2}=q$ and
$q_{1,2}q_{2,1}=q^{-1}$, with a primitive root of unity $q$.  From the
logarithmic perspective, this $W_3$ algebra is a \textit{nonextended}
algebra, playing the same role in relation to an extended algebra as
the Virasoro algebra plays in relation to the triplet algebras of
$(p,1)$ \cite{[Kausch],[Gaberdiel-K],[FGST]} and $(p,p')$
\cite{[FGST3]} logarithmic models.  Specifically in the $W_3$ case,
the \textit{extended} algebra is the octuplet algebra described in
Appendix~\ref{app:W}.  Similar constructions are expected in other
good cases; I am optimistic about the fact that the same generalized
Dynkin diagram gives rise to a finite-dimensional Nichols algebra
\textit{and} to an interesting conformal field theory.  The intricate
machinery underlying the finite dimensionality of the corresponding
Nichols algebra may manifest itself in constructing new logarithmic
models.\footnote{Recall that \textit{rational} conformal field
  theories are generally defined as the cohomology of a complex
  associated with the screenings, whereas \textit{logarithmic} models
  are defined by the kernel (cf.~\cite{[FHST],[FGST],[FGST3],[AM-3]}).
  In particular, this allows interesting logarithmic conformal models
  to exist in the cases where the rational model is nonexistent (the
  $(p,1)$ series) or trivial (the $(2,3)$ model).}



In what follows, I therefore reproduce Heckenberger's list of rank-2
finite-dimensional Nichols algebras~\cite{[Heck-1+1]}, with the only
difference that I \texttt{enumerate}, not \texttt{itemize} the
subitems.  For several items, I also add the presentations known
from~\cite{[Ag-0804-standard]} and explicitly borrowed
from~\cite{[Hel]}, including the case number in~\cite{[Hel]}.  From
the Nichols-algebra data, I move toward conformal field theory by
analyzing the conditions on the screening momenta.  When it is clear
what current algebra extends the Virasoro algebra with the central
charge obtained from~\eqref{c-rank2}, I recall the explicit
construction, presenting it in the form that manifestly refers to the
corresponding pair of screenings (once again, all extended algebras
except the one in Appendix~\ref{app:W} are not logarithmic extensions,
but rather starting points for such extensions).

\subsection*{0.9. \ Points to note}
\begin{enumerate}\renewcommand{\labelenumi}{0.9\theenumi.}%
  \addtolength{\itemsep}{6pt}%
\item In conformal field theory, \textit{fermionic} screenings are
  often interesting.
  Their Nichols-algebra counterparts are the diagonal entries $-1$ in
  braiding matrices.  But given a $q_{i,i}=-1$ and trying to
  reconstruct a screening in general leads to the condition
  $\alpha_i.\alpha_i=1+2m$ on the screening momentum, with $m\in\oZ$.
  For the corresponding screening current
  $f(z)=e^{\alpha_i.\varphi(z)}$, it then follows that $f(z)f(w)$
  develops a $(1+2m)$th-order zero as $z\to w$.  The cases where this
  zero is actually a pole are somewhat pathological from the CFT{}
  standpoint; as regards the cases of a zero of an order $\geq 3$, I
  am unaware of any such examples of screenings.  Only $m=0$ is a
  ``good'' value.  Remarkably, \textit{solving
    conditions~\eqref{cartan-test-log} with $q_{i,i}=-1$ has the
    tendency to select the value $m=0$, thus ensuring a true fermionic
    screening}.

\item Other integers appearing in ``taking the logarithms'' are not
  disposed of that easily.  There are solutions
  of~\eqref{cartan-test-log} where these integers vanish (and the
  central charge depends on another integer parameter, the
  \textit{order} of a root of unity); such solutions are referred to
  as ``regular'' in what follows.  But there also exist ``peculiar''
  solutions of~\eqref{cartan-test-log} where some of these parasitic
  integers persist, and which have somewhat reduced chances to
  correspond to interesting CFT{} models.  In fact, \textit{some}
  peculiar solutions are eliminated already by the conditions in
  Heckenberger's list: in some items, the order of the corresponding
  root of unity must not be too small, and the peculiar solutions do
  require just one of those excluded values.  This might suggest that
  peculiar solutions should somehow be eliminated altogether, but if
  so, then I have overlooked the argument.

\item Things get worse with the many items in the list that do not
  involve a free discrete parameter such as the order of a root of
  unity.  Isolated central charge values are by no means illuminating,
  and remain entirely unsuggestive when expanded into families by the
  occurrence of ``parasitic integers.''

  The unwieldiness of the ``peculiar'' central charges also thwarted
  my original intention to provide each item in the list with a
  central charge.  This can be done, but the results are not
  indicative of anything.  The corresponding items in the list are
  therefore left in their original form given in~\cite{[Heck-1+1]}.

\item In the ``regular'' cases, I choose a primitive $p$th root of
  unity as $e^{\frac{2i\pi}{p}}$ with an integer $p$.
  This might unnecessarily restrict the generality, but the cases that
  follow with this choice are already interesting.  In ``peculiar''
  cases, by contrast, I try to work out the cases with $e^{\frac{2i\pi
      r}{p}}$, where $r$ is coprime with $p$.  The $r$ parameter
  sometimes survives till the central charge, but that's where the
  story ends, because I do not construct any current algebra
  generators beyond Virasoro in peculiar cases.

  Mostly, I take the logarithm of relations such as $e^{i\pi
    x}=e^{\frac{2i\pi r}{s}}$ (where $x$ is typically a linear
  combination of scalar products) ``honestly,'' as $x=\frac{2r}{s}+2
  \ell$, $\ell\in\oZ$.  In \textit{some} cases, the ensuing dependence
  on $\ell$ turns out to be ``under control'' (something like a shift
  of the level of an affine Lie algebra with which the corresponding
  conformal field theory is associated---which interestingly
  corresponds to a \textit{twist} equivalence of the braiding matrix),
  and I sometimes omit it.

\item Strictly speaking, identifying a CFT{} model from its central
  charge that depends on a parameter is an ill-defined procedure in
  the sense that given a central charge $c=f(p)$ and redefining the
  parameterization by an arbitrary function, $p'=g(p)$, changes the
  ``functional form'' of $c$ arbitrarily.  It is tacitly understood
  that some ``natural'' parameterizations are considered and very
  limited reparameterizations are allowed (typically those that are
  known to occur in some CFT{} constructions).
\end{enumerate}

\subsection*{0.99. \ Notation.}
The notation 
\begin{equation*}
  R_{\ell} ={}\text{the set of primitive $\ell$th roots of unity}
\end{equation*}
is copied from \cite{[Heck-1+1]} as part of the defining conditions in
the list items.  A braiding matrix~\eqref{q-matrix} is encoded in a
generalized Dynkin diagram $\Dynkin{q_{1,1}}{m_{1,2}}{\ \
  q_{2,2}}$\!\!\!, where $m_{1,2}=q_{1,2} q_{2,1}$.  The $A_2$, $B_2$,
and $G_2$ Cartan matrices are $\left(\begin{smallmatrix}
    2&-1\\
    -1&2
  \end{smallmatrix}\right)$, $\left(\begin{smallmatrix}
    2&-2\\
    -1&2
  \end{smallmatrix}\right)$, and $\left(\begin{smallmatrix}
    2&-3\\
    -1&2
  \end{smallmatrix}\right)$.
Vectors $\alpha,\beta\in\oC^2$ are the momenta of two screenings in
what follows ($\alpha_1$ and $\alpha_2$ in the nomenclature
of~\eqref{FF}).

\section{The list, item \textbf{1}}
The defining conditions are
\begin{equation*}
  \boxed{q_{12}q_{21}=1
  \quad\text{and}\quad
  q_{11},q_{22}\in \bigcup_{a=2}^\infty R_a.}
\end{equation*}
This is the ``trivial'' $A_1\times A_1$ case.  The corresponding CFT{}
model is the product $(p',1)\times(1,p))$ of two ``$(p,1)$''
models~\cite{[FGST],[FGST2]}, or, in the degenerate case where
$\alpha$ and $\beta$ are collinear (and hence only one boson is
needed), the $(p',p)$ model~\cite{[FGST3],[FGST4]}.

\section{The list, items \textbf{2.}$\pmb{*}$}%
\renewcommand{\arraycolsep}{0pt}%
The defining conditions are
\begin{equation*}
  \boxed{q_{12}q_{21}q_{22}=1 \quad\text{and}\quad q_{12}q_{21}\neq1,}
\end{equation*}
plus any of conditions \textbf{2.1}--\textbf{2.7}.  In terms of the
momenta $\alpha,\beta\in\oC^2$ of the screenings, the common condition
for all these cases takes the form
\begin{equation*}
  2\alpha.\beta + \beta.\beta=2m\ \ (m\in\oZ).
\end{equation*}

\subsection*{2.1 \ (\Hel{5.7}{1}).} $q_{11}q_{12}q_{21}=1$, \
$q_{12}q_{21}\in \bigcup_{a=2}^\infty R_a$, \ Cartan type $A_2$,
$\Dynkin{q}{q^{-1}}{q}$.

In terms of scalar products, the conditions are
\begin{equation*}
  \alpha.\alpha + 2\alpha.\beta=2n\ \ (n\in\oZ),
  \quad 2\alpha.\beta=-\ffrac{2}{p} + 2 j,\quad |p|\geq 2\ \ (j\in\oZ).
\end{equation*}\pagebreak[3]%
The braiding matrix (which is stable under Weyl reflections) is then
parameterized as {\footnotesize$\begin{pmatrix}
    e^{\frac{2 i \pi }{p}} & (-1)^j e^{-\frac{i \pi}{p}} \\
    (-1)^j e^{-\frac{i \pi}{p}} & e^{\frac{2 i \pi }{p}}
  \end{pmatrix}$}. \ None of the two screenings is fermionic unless
$|p|=2$.

Conditions~\eqref{cartan-test-log} are not satisfied for all $m$, $n$,
$j$, and $p$.  There are several ``peculiar'' solutions and a
``regular'' solution.  The peculiar solutions are $(m =
0,\linebreak[0]\,n = -k-\frac{3}{2},\linebreak[0]\,p =
2,\linebreak[0]\,j = -k-\frac{3}{2})$, \ $(m = 0,\linebreak[0]\,n =
-k-\frac{3}{2},\linebreak[0]\,p = -2,\linebreak[0]\,j =
-k-\frac{5}{2})$, \ $(m = -k-\frac{3}{2},\linebreak[0]\,n =
0,\linebreak[0]\,p = 2,\linebreak[0]\,j = -k-\frac{3}{2})$, \ $(m =
-k-\frac{3}{2},\linebreak[0]\,n = 0,\linebreak[0]\,p =
-2,\linebreak[0]\,j = -k-\frac{5}{2})$, \ $(m =
k+\frac{3}{2},\linebreak[0]\,n = k+\frac{3}{2},\linebreak[0]\,p =
2,\linebreak[0]\,j = k+\frac{3}{2})$, and $(n = m =
k+\frac{3}{2},\linebreak[0]\, k+\frac{3}{2},\linebreak[0]\,p =
-2,\linebreak[0]\,j = k+\frac{1}{2})$ with half-integer~$k$ in all
cases; with this parameterization, the resulting central charge is in
each case equal to $\fffrac{3k}{k+2}-1$, which is the central charge
of the $\hSL2_k/\heis$ coset (more on it is to be said below, when it
occurs as a ``regular'' solution).

The regular solution is $m=n=0$, yielding the central charge
\begin{equation}\label{W3-cc}
  c= 50 - \ffrac{24}{k+3} - 24 (k+3),
\end{equation}
where $k+3 = \frac{1}{p} -j$ (or, in view of the structure of the
formula, $\frac{1}{k+3} = \frac{1}{p} - j$).  This is the central
charge of the $W_3$ algebra parameterized in terms of the level $k$ of
the $\hSL3$ affine Lie algebra from which $W_3$ can be obtained by
Hamiltonian reduction.

The centralizer of the screenings does indeed contain a dimension-3
primary field $W(z)$ (unique up to an overall factor) in the space of
differential polynomials in the fields
\begin{equation*}
  \partial\varphi_\alpha(z)=\alpha.\partial\varphi(z)
  \quad\text{and}\quad
  \partial\varphi_\beta(z)=\beta.\partial\varphi(z).
\end{equation*}
Explicitly, setting $j=0$ for simplicity and omitting the $(z)$
arguments in the right-hand side for brevity,
\begin{multline}\label{W3-gen}
  W(z) =
  \partial\varphi_\alpha \partial\varphi_\alpha \partial\varphi_\alpha
  +
  \ffrac{3}{2}\partial\varphi_\alpha \partial\varphi_\alpha
  \partial\varphi_\beta 
  -\ffrac{3}{2} \partial\varphi_\alpha \partial\varphi_\beta
  \partial\varphi_\beta
  -\partial\varphi_\beta \partial\varphi_\beta \partial\varphi_\beta
  \\
  -\ffrac{9 (p - 1)}{2
    p}\partial^2\varphi_\alpha \partial\varphi_\alpha -\ffrac{9 (p -
    1)}{4 p} \partial^2\varphi_\alpha \partial\varphi_\beta + \ffrac{9
    (p - 1)}{4 p}\partial^2\varphi_\beta \partial\varphi_\alpha +
  \ffrac{9 (p - 1)}{2 p} \partial^2\varphi_\beta \partial\varphi_\beta
  \\
  + \ffrac{9 (p - 1)^2}{4 p^2} \partial^3\varphi_\alpha -\ffrac{9 (p -
    1)^2}{4 p^2} \partial^3\varphi_\beta.
\end{multline}

The Nichols algebra $\BX$ (of the two-dimensional braided vector space
$X$ with basis $F_1$ and $F_2$ with braiding matrix~\eqref{q-matrix})
is in this case the quotient~\cite{[Hel]}
\begin{equation}\label{Nich-W3}
  \BX=T(X)/ \bigl([F_1,F_1,F_2],\ [F_1,F_2,F_2],\ F_1^{p},\
  [F_1,F_2]^{p},\ F_2^{p}\bigr)
\end{equation}
if $p\geq 3$.  Here and hereafter, square brackets denote iterated
$q$-commutators constructed in accordance with Lyndon word
decomposition.  If $p=2$ (the screenings \textit{are} fermionic!), the
triple-bracket generators of the ideal are absent.  The dimension is
$\dim\Nich(X)=p^3$.

The elements $F_1^{p}$ and $F_2^{p}$ in the ideal indicate the
``positions'' of the \textit{long screenings}
\begin{equation}\label{long-scr}
  \Ea=\oint e^{-\alpha^{\vee}.\varphi}
  =\oint e^{-p\alpha.\varphi}
  \quad\text{and}\quad
  \Eb=\oint e^{-\beta^{\vee}.\varphi}
  =\oint e^{-p\beta.\varphi}
  \qquad
  (\alpha^{\vee}=\fffrac{2\alpha}{\alpha.\alpha}),
\end{equation}
i.e., $F_1^{p}$ and $F_2^{p}$ ``tend to be'' the operators
``opposite'' to the respective long screening.  Generally, these long
screenings are to produce $m$-plet structures in logarithmic models,
similarly to how the triplet structure of the $(p,1)$ logarithmic
models~\cite{[Kausch],[Gaberdiel-K]} is generated by the corresponding
long screening~\cite{[FHST]}.  For the current $W_3$-case, the
resulting \textit{octuplet} algebra is outlined in
Appendix~\ref{app:W}.

\subsection*{2.2\ \ (\Hel{5.7}{3}).} $q_{11}=-1$, \ $q_{12}q_{21}\in
\bigcup_{a=3}^\infty R_a$, \ Cartan type $A_2$, \
$\Dynkin{-1}{q^{-1}}{q}$.

The conditions are restated in terms of the scalar products as
\begin{equation*}
  \alpha.\alpha = 1 + 2 n\ \ (n\in\oZ),
  \quad 2\alpha.\beta=-\ffrac{2}{p} + 2 j, \ \ |p|\geq 3\ \ (j\in\oZ).
\end{equation*}
The braiding matrix is then parameterized as
{\footnotesize$\begin{pmatrix}
    -1 & (-1)^j e^{- \frac{i \pi}{p}} \\
    (-1)^j e^{- \frac{i \pi}{p}} & e^{\frac{2 i \pi }{p}}
  \end{pmatrix}$}, and its Weyl reflections different from the
original matrix is {\footnotesize$\begin{pmatrix}
    -1 & -(-1)^j e^{\frac{i\pi}{p}} \\
    -(-1)^j e^{\frac{i\pi}{p}} & -1
\end{pmatrix}$}.
The Weyl orbit also contains the braiding matrix
{\footnotesize$\begin{pmatrix}
    e^{\frac{2i\pi}{p}} & (-1)^j e^{-\frac{i\pi}{p}} \\
 (-1)^j e^{  -\frac{i\pi}{p}} & -1
\end{pmatrix}$}.

The first screening ``wants to be fermionic.''  Remarkably,
conditions~\eqref{cartan-test-log} have a solution only if $m=n=0$
(would-be solutions with nonzero $m$ or $n$ require $|p|=2$).  Thus
the $-1$ in the braiding matrix does indeed correspond to a fermionic
screening in the standard sense---an operator of the form $\oint
F(z)$, where $F(z)F(w)$ has a first-order, not a higher-order, zero as
$z\to w$.

The solution for the scalar products with $m=n=0$ yields the central
charge
\begin{equation*}
  c= \ffrac{3 k}{k + 2} - 1,
\end{equation*}
where $k + 2 = \frac{1}{p} -j$.  This is the central charge of the
$\hSL2_k/\heis$ coset (where $\heis$ is the Heisenberg subalgebra).
The two currents $\jplus(z)$ and $\jminus(z)$ that are in the
centralizer of the screenings and generate the coset algebra can be
expressed in terms of the two bosons ``in the direction'' of each
screening as\footnote{The exponentials in~\eqref{sl2:nonsymm} are
  assumed to be normal ordered, and the second line involves the
  (standard) abuse of notation: nested normal ordering from right to
  left is in fact understood after the expression is expanded.}
\begin{equation}\label{sl2:nonsymm}
  \begin{aligned}
    \jplus(z)& = e^{-\frac{1}{k}(2\varphi_{\alpha}(z) +
      \varphi_{\beta}(z))},\\
    \jminus(z)& =
    -(\partial\varphi_{\alpha}(z)\partial\varphi_{\beta}(z)
    +\partial\varphi_{\alpha}(z)\partial\varphi_{\alpha}(z)
    +(k+1)\partial^2\varphi_{\alpha}(z))
    e^{\frac{1}{k}(2\varphi_{\alpha}(z) + \varphi_{\beta}(z))}.
  \end{aligned}
\end{equation}
Adding a boson $\chi(z)$ associated with the $\heis$ algebra
immediately yields the three $\hSL2_k$ currents
\begin{equation*}
  J^{\pm}(z) = j^{\pm}(z) e^{\pm\sqrt{\frac{2}{k}}\chi(z)}
  \quad\text{and}\quad
  \Jnaught(z) = \sqrt{\ffrac{k}{2}}\,\partial\chi(z).
\end{equation*}

For $p\geq 3$, the Nichols algebra is the quotient~\cite{[Hel]}
\begin{equation*}
  \BX=T(X)/ \bigl([F_1,F_2,F_2],\ F_1^2,\ F_2^{p}\bigr),
\end{equation*}
with $\dim\BX=4p$.

A long screening here is
\begin{equation*}
  \mathscr{E}_\beta = \oint e^{-p\beta.\varphi}.
\end{equation*}

\subsection*{2.3\ (\Hel{5.11}{3})} $q_{11}\in R_3$, \ $q_{12}q_{21}\in
\bigcup_{a=2}^\infty R_a$, \ $q_{11}q_{12}q_{21}\neq1$, \ Cartan type
$B_2$, \ $\Dynkin{\zeta}{q^{-1}}{q}$, \ $R_3\ni\zeta\neq q$.

In terms of the screening momenta, the conditions become
\begin{equation*}
  \alpha.\alpha = \ffrac{2s}{3},
  \quad
  2\alpha.\beta = -\ffrac{2}{p} + 2 j,\ \ |p|\geq 2\ \ (j\in\oZ),
\end{equation*}
where $s$ is coprime with~$3$.  The braiding matrix is parameterized
as {\footnotesize$\begin{pmatrix}
    e^{\frac{2 i \pi  s}{3}} & (-1)^j e^{-\frac{i \pi}{p}} \\
    (-1)^j e^{-\frac{i \pi}{p}} & e^{\frac{2 i \pi }{p}}
  \end{pmatrix}$} and its Weyl reflection noncoincident with the
original is {\footnotesize$\begin{pmatrix}
    e^{\frac{2 i \pi  s}{3}} & (-1)^j e^{-\frac{i \pi  (4 p s-3)}{3 p}} \\
    (-1)^j e^{-\frac{i \pi (4 p s-3)}{3 p}} & e^{\frac{2 i \pi (4 p
        s-3)}{3 p}}
  \end{pmatrix}$}.  Conditions~\eqref{cartan-test-log} can be
satisfied only if $(m = 0, p = 3, s = 3 \ell - 1, j = 1 - 2 \ell)$ or
$(m = 0, p = -3, s = 3 \ell + 1, j = -1 - 2 \ell)$ (two peculiar
solutions), or $(m = 0, s = 1)$ (the regular solution).

In both peculiar cases, the central charge is
  $86 - 60 (k+3) - \ffrac{30}{k+3}$,
where $k+3 = -\frac{1}{3} + \ell$ (or $\frac{1}{k+3}=-\frac{2}{3} + 2
\ell$) in the first case and $k+3=\frac{1}{3} + \ell$ (or
$\frac{1}{k+3}=\frac{2}{3} + 2 \ell$) in the second case.  The central
charge is that of the $WB_2$ algebra, discussed in more detail below
when it appears as a ``regular'' solution.

In the regular case $m=0$, the central charge is
\begin{equation*}
  c=-\ffrac{21}{2} - 6 z - \ffrac{27}{2 (4 z -3)},
\end{equation*}
where $\frac{1}{z}=\frac{1}{p} -j$ or $\frac{1}{z}= -\frac{1}{p} +
\frac{4}{3} + j$.

The condition $q_{11}q_{12}q_{21}\neq1$
excludes the value $p=3$.

If $p\geq 4$,
and
$p'=\ord(q_{11}q_{22}^{-1})=\ord(e^{\frac{2i\pi}{3}-\frac{2i\pi}{p}})$,
then~\cite{[Hel]}
\begin{equation*}
  \BX=T(X)\!\bigm/\!\bigl( [F_1,F_2,F_2],\ F_1^{3},\ [F_1,F_1,F_2]^{p'},\
  F_2^{p}\bigr),
\end{equation*}
with $\dim\BX=9 p p'$.

\subsection*{2.4} $q_{11}\in \bigcup _{n=4}^\infty R_a$, with two
subcases listed below.  To identify the central charges in what
follows, we use the formula~\cite{[FKW]}
\begin{equation}\label{FKW}
  c(k) = \ell - 12\frac{|(k+h^{\vee})\rho^{\vee} - \rho|^2}{k+h^{\vee}}
\end{equation}
for the central charge of a $W$-algebra obtained by Hamiltonian
reduction of a level-$k$ affine Lie algebra; $h^{\vee}$ is the dual
Coxeter number, $\rho$ is half the sum of positive roots,
$\rho^{\vee}$ half the sum of their duals, and $\ell$ is the rank of
the corresponding finite-dimensional Lie algebra.

\begin{description}
\item[2.4.1\ \ (\Hel{5.11}{1})] $q_{12}q_{21} = q_{11}^{-2}$, \ Cartan
  type $B_2$, \ $\Dynkin{q}{q^{-2}}{q^2}$\!\!.
  
  In terms of scalar products, we then have
  \begin{equation*}
    \alpha.\alpha=\ffrac{2}{p} + 2 j,\ \ |p|\geq4\ \ (j\in\oZ),\quad
    2\alpha.\beta + 2\alpha.\alpha = 2n\ \ (n\in\oZ).
  \end{equation*}
  The braiding matrix (stable under Weyl reflections) is
  {\footnotesize$\begin{pmatrix}
      e^{\frac{2 i \pi }{p}} & (-1)^n e^{-\frac{2i \pi}{p}} \\
      (-1)^n e^{-\frac{2 i \pi}{p}} & e^{\frac{4 i \pi }{p}}
    \end{pmatrix}$}.  Conditions~\eqref{cartan-test-log} hold in two
  peculiar cases and one regular case.  The peculiar cases are $(m =
  -2 j, n = 0, p = 4)$ with $c=-1-\frac{24}{4 j+1}+\frac{24}{4 j-1}$
  and $(m = 1 - 2 j, n = 0, p = -4)$ with $c=-1-\frac{24}{4 j-1}
  +\frac{24}{4 j-3}$.

  The regular case is $m=n=0$, with
  \begin{equation*}
    c = 86 - 60 (k+3) - \ffrac{30}{k+3},
  \end{equation*}
  where $k+3 =\frac{1}{p} + j$ (or $\frac{1}{k+3} = \frac{2}{p} + 2
  j$).  This is the central charge of the $WB_2$
  algebra~\cite{[spin-4],[FfSchTh],[IT]} (also see~\cite{[KW]})
  obtained by Hamiltonian reduction of the level-$k$ \ $B^{(1)}_2$
  algebra (by formula~\eqref{FKW}, with $|\rho|^2=\frac{5}{2}$,
  $|\rho^{\vee}|^2=5$, and
  $\langle\rho,\rho^{\vee}\rangle=\frac{7}{2}$ for $B_2$).

  The $WB_2$ algebra contains a unique primary field of dimension~$4$.
  Explicitly, it is a rather long (20 terms) differential polynomial
  in $\partial\varphi_{\alpha}(z)$ and $\partial\varphi_{\beta}(z)$,
  \begin{multline*}
    p (p-3) (27 p-32) \partial\varphi_{\alpha} \partial\varphi_{\alpha}
    \partial\varphi_{\alpha} \partial\varphi_{\alpha}
    + 2 p (p-3) (27 p-32) \partial\varphi_{\alpha} \partial\varphi_{\alpha}
    \partial\varphi_{\alpha} \partial\varphi_{\beta}\\
    {}-21 p
    \left(p^2-2\right) \partial\varphi_{\alpha} \partial\varphi_{\alpha}
    \partial\varphi_{\beta} \partial\varphi_{\beta}
    - p (3 p-2) (16
    p-27) \partial\varphi_{\alpha} \partial\varphi_{\beta}
    \partial\varphi_{\beta} \partial\varphi_{\beta} \\
    \shoveright{{}-\fffrac{p}{4} (3 p-2) (16
      p-27) \partial\varphi_{\beta} \partial\varphi_{\beta}
      \partial\varphi_{\beta} \partial\varphi_{\beta}}
    \\
    +\dots\dots\dots\dots\dots\dots\dots\dots\dots\dots\dots\\
    - \fffrac{(p\!-\!3) (3 p\!-4\!) (30 p^3\!-\!115 p^2\!+\!144
      p\!-\!60)}{3 p^2} \partial^4\varphi_{\alpha} +\fffrac{(2
      p\!-\!3) (3 p\!-\!2) (15 p^3\!-\!72 p^2\!+\!115 p\!-\!  60)}{3
      p^2}\partial^4\varphi_{\beta}
  \end{multline*}
  (\textit{all} coefficients are polynomials in $p$ with integer
  coefficients after the overall renormalization by $12p^2$).

  If $p\geq 5$ is odd, then~\cite{[Hel]}
  \begin{equation*}
    \BX=T(X)\!\bigm/\!\bigl( [F_1,F_1,F_1,F_2] , \ [F_1,F_2,F_2] , \
    F_1^{p} ,\
    [F_1,F_1,F_2]^{p},\ [F_1,F_2]^{p} , \ F_2^{p} \bigr),
  \end{equation*}
  with $\dim\BX=p^4$.  If $p\ge 4$ is even, then
  \begin{equation*}
    \BX=T(X)\!\bigm/\!\bigl( [F_1,F_1,F_1,F_2] , \ [F_1,F_2,F_2] , \
    F_1^{p} ,\
    [F_1,F_1,F_2]^{\frac{p}{2}},\ [F_1,F_2]^{p} , \ F_2^{\frac{p}{2}}
    \bigr)
  \end{equation*}
  and $\dim\BX=\frac{p^4}{4}$ (the second generator of the ideal is
  absent for $p=4$).

\medskip

\item[2.4.2] $q_{12}q_{21} = q_{11}^{-3}$, \ Cartan type $G_2$, \
  $\Dynkin{q}{q^{-3}}{q^3}$.

  In terms of scalar products, we now have
  \begin{equation*}
    \alpha.\alpha=\ffrac{2}{p} + 2 j,\ \ |p|\geq4\ \ (j\in\oZ),\quad
    2\alpha.\beta + 3\alpha.\alpha = 2n\ \ (n\in\oZ).
  \end{equation*}
  The braiding matrix is parameterized as {\footnotesize$\begin{pmatrix}
      e^{\frac{2 i \pi }{p}} & (-1)^{j+n} e^{ -\frac{3i \pi}{p}} \\
      (-1)^{j+n} e^{ -\frac{3i \pi}{p}} & e^{\frac{6 i \pi }{p}}
    \end{pmatrix}$} and is stable under Weyl reflections.
  Conditions~\eqref{cartan-test-log} hold in the peculiar cases $(j =
  0, m = 0, p = 4)$ with $c=-10-\frac{54}{4 n+1}+\frac{24}{4 n-3}$,
  $(m = -3 j, n = 0, p = 6)$ with $c=-\frac{2}{3} + \frac{400}{3 (18
    j-1)} -\frac{36}{6 j+1}$, and $(m = 1 - 3 j, n = 0, p = -6)$ with
  $c=-\frac{2}{3} + \frac{400}{3 (18 j-7)}-\frac{36}{6 j-1}$ and in
  the regular case $m = n = 0$ with the central
  charge
  \begin{equation*}
    c = 194 - 168(k+4) - \ffrac{56}{k+4},
  \end{equation*}
  where $k+4=\frac{1}{p} + j$ (or $\frac{1}{k+4} = \frac{3}{p} + 3
  j$).  This is the central charge of the $WG_2$
  algebra~\cite{[spin-6],[IT]} (also see~\cite{[KW],[Zhu]}) obtained
  by Hamiltonian reduction of the level-$k$ \ $G^{(1)}_2$ algebra (by
  formula~\eqref{FKW}, with $|\rho|^2=14$,
  $|\rho^{\vee}|^2=\frac{14}{3}$, and
  $\langle\rho,\rho^{\vee}\rangle=8$ for~$G_2$).  The $WG_2$ algebra
  contains a unique primary field of dimension~$6$, which is by far
  too long to be given here (see~\cite{[IT],[spin-6],[Zhu]}).
\end{description}

\medskip

The remaining subcases of Case~2 may all be considered ``peculiar'' to
some extent.  The values of $c$ are equally ``peculiar.''

\subsection*{2.5} $q_{12}q_{21}\in R_8$, \ $q_{11}=(q_{12}q_{21})^2$,
\ Cartan type $G_2$, \ $\Dynkin{\zeta^2}{\zeta}{\zeta^{-1}}$\!\!\!, \
$\zeta\in R_8$.

The conditions reformulate in terms of scalar products as
\begin{equation*}
  2\alpha.\beta = \ffrac{r}{4} + 2 j\ \ (j\in\oZ),\quad
  \alpha.\alpha - 4\alpha.\beta = 2n,
\end{equation*}
where $r$ is odd.  The braiding matrix is {\footnotesize$\begin{pmatrix}
  i^r & (-1)^j e^{\frac{i \pi  r}{8}} \\
  (-1)^j e^{\frac{i \pi r}{8}} & e^{-\frac{1}{4} i \pi r}
\end{pmatrix}$} and its Weyl reflection different from the original
matrix is
{\footnotesize$\begin{pmatrix} i^r & (-1)^j e^{\frac{3}{8} i \pi r} \\
  (-1)^j e^{\frac{3}{8} i \pi r} & (-1)^r
\end{pmatrix}$}.
Conditions~\eqref{cartan-test-log} can be satisfied only if $(m = 0, r
= 1 - 8 j - 4 n)$, yielding the central charge $c=-10 - \frac{48}{4
  n-1} + \frac{108}{4 n-9}$.  For $n=0$, this gives the celebrated
central charge value
\begin{equation*}
  c=26.
\end{equation*}

\subsection*{2.6} $q_{12}q_{21}\in R_{24}$, \
$q_{11}=(q_{12}q_{21})^6$, \
$\Dynkin{\zeta^6}{\zeta}{\zeta^{-1}}$\kern-5pt, \ $\zeta\in R_{24}$.

The conditions for the scalar products are
\begin{equation*}
  2\alpha.\beta = \ffrac{r}{12} + 2 j\ \ (j\in\oZ),\quad
  \alpha.\alpha - 12\alpha.\beta = 2n,
\end{equation*}
where $r$ is coprime with $2$ and $3$.  The braiding matrix
{\footnotesize$\begin{pmatrix}
    i^r & (-1)^j e^{\frac{i \pi  r}{24}} \\
    (-1)^j e^{\frac{i \pi r}{24}} & e^{-\frac{i \pi r}{12}}
  \end{pmatrix}$} has the $G_3$ Cartan matrix, but its nontrivial Weyl
reflection is {\footnotesize$\begin{pmatrix} i^r & (-1)^j e^{\frac{11}{24} i \pi r} \\
    (-1)^j e^{\frac{11}{24} i \pi r} & e^{\frac{2 i \pi r}{3}}
\end{pmatrix}$}, with the associated generalized Cartan matrix
$\left(\begin{smallmatrix} 2& -3\\ -2& 2
  \end{smallmatrix}\right)$; various other generalized Cartan
matrices are produced under further Weyl reflections.

Conditions~\eqref{cartan-test-log} are satisfied only if $(m = 0, r =
1 - 24 j - 4 n)$, with $c=-10 - \frac{144}{4 n-1}+\frac{324}{4 n-25}$
(that $r$ be coprime with $2$ and $3$ selects the values $n=2 + 3
\ell$ or $n=3 + 3 \ell$, $\ell\in\oZ$).

\subsection*{2.7} $q_{12}q_{21}\in R_{30}$, \
$q_{11}=(q_{12}q_{21})^{12}$,
$\Dynkin{\zeta^{12}}{\zeta}{\zeta^{-1}}$\!\!\!,\ \ $\zeta\in R_{30}$.

In terms of scalar products,
\begin{equation*}
  2\alpha.\beta = \ffrac{r}{15} + 2 j\ \ (j\in\oZ),\quad
  \alpha.\alpha -24 \alpha.\beta = 2n,
\end{equation*}
where $r$ is coprime with $30$.  The braiding matrix, parameterized
as
{\footnotesize$\begin{pmatrix} e^{\frac{4 i \pi r}{5}} & (-1)^j e^{\frac{i \pi  r}{30}} \\
    (-1)^j e^{\frac{i \pi r}{30}} & e^{-\frac{i \pi r}{15}}
  \end{pmatrix}$}, has the associated generalized Cartan matrix
$\left(\begin{smallmatrix}
    2 & -4 \\
    -1 & 2
  \end{smallmatrix}\right)$.  The nontrivial Weyl reflection is
{\footnotesize$\begin{pmatrix} e^{\frac{4 i \pi r}{5}} &
    (-1)^{j+r} e^{-\frac{7}{30} i \pi  r} \\
    (-1)^{j+r} e^{-\frac{7}{30} i \pi r} & (-1)^r
  \end{pmatrix}$}, with the same generalized Cartan matrix, but other
(generalized) Cartan matrices are produced under further Weyl
reflections.  Conditions~\eqref{cartan-test-log} are solved by $(m =
0, n = 1 + 2 \ell, r = -2 - 5 \ell - 30 j)$, where $\ell = 1 + 6 u$ or
$\ell = 3 + 6 u$, $u\in\oZ$, respectively yielding the
incomprehensible $c=-\frac{62}{5} + \frac{2916}{5 (30 u -17)} -
\frac{180}{30 u + 7}$ and $c=-\frac{62}{5} + \frac{2916}{5 (30 u-7)} -
\frac{180}{30 u + 17}$.


\section{The list, items \textbf{3.}$\pmb{*}$}
The defining conditions are
\begin{equation*}
  \boxed{q_{12}q_{21}\neq1, \ \ q_{11}q_{12}q_{21}\neq1, \ \
    q_{12}q_{21}q_{22}\neq1, \ \ q_{22}=-1, \ \ q_{11}\in R_2\cup R_3,}
\end{equation*}
plus any of conditions \textbf{3.1}--\textbf{3.7}.  In terms of the
screening momenta, the common conditions are
\begin{equation*}
  \beta.\beta=1+2m \ \ (m\in\oZ),\quad
  \alpha.\alpha=1\text{ or }\ffrac{2 s}{3},
\end{equation*}
where $s$ is coprime with~$3$.

\subsection*{3.1 \ (\Hel{5.7}{4})} $q_{11}=-1$, \ $q_{12}q_{21}\in
\bigcup_{a=3}^\infty R_a$, \ Cartan type $A_2$, \
$\Dynkin{-1}{q}{-1}$.

In terms of scalar products of the screening momenta, these conditions
are
\begin{equation*}
  \alpha.\alpha=1 + 2 n \ \ (n\in\oZ), \quad
  2\alpha.\beta=\ffrac{2}{p} + 2 j,\ \ |p|\geq 3\ \ (j\in\oZ).
\end{equation*}
Both screenings ``want to be fermionic.''  The braiding matrix is
{\footnotesize$\begin{pmatrix}
    -1 & (-1)^j e^{\frac{i \pi }{p}} \\
    (-1)^j e^{\frac{i \pi }{p}} & -1
  \end{pmatrix}$} and its Weyl reflections are {\footnotesize$\begin{pmatrix}
    -1 & -(-1)^j e^{-\frac{i \pi}{p}} \\
    -(-1)^j e^{-\frac{i \pi}{p}} & e^{\frac{2 i \pi }{p}}
   \end{pmatrix}$} and {\footnotesize$\begin{pmatrix} e^{\frac{2 i \pi }{p}}
     & -(-1)^j e^{-\frac{i \pi}{p}}
     \\
     -(-1)^j e^{-\frac{i \pi}{p}} & -1
\end{pmatrix}$}.

Remarkably, conditions~\eqref{cartan-test-log} are
satisfied (with $|p|\geq3$) only for $m=n=0$ (no ``peculiar''
solutions!), yielding the $\hSL2_k/\heis$ central charge
\begin{equation*}
  c= \ffrac{3 k}{k + 2} - 1,
\end{equation*}
where $k + 1 = \frac{1}{p} + j$.  For $j=0$, in particular, this
relation between $k$ and $p$ takes the form
\begin{equation*}
  \mfrac{1}{p + 1} + \mfrac{1}{k + 2} = 1.
\end{equation*}
This ``duality'' between two levels, $k$ and $p-1$, was extensively
used in~\cite{[BFST],[W2n]} (see also the references therein); in
particular,
\begin{equation*}
  \hSL2_k/\heis=\widehat{s\ell}(2|1)_{p-1}/\widehat{g\ell}(2)_{p-1},
\end{equation*}
offering another view on what the CFT{} counterpart of the Nichols
algebra is.\footnote{This coset equivalence belongs to a vast subject
  discussed in~\cite{[W-uni]}.}

The currents generating the $\hSL2_k/\heis$ coset algebra are given by
\begin{equation}
  \label{sl2:symm}
  \begin{aligned}
    \jplus(z) &= \partial\varphi_{\beta}(z)\,
    e^{\frac{1}{k}(\varphi_{\alpha}(z)-\varphi_{\beta}(z))},\\
    \jminus(z) &= \partial\varphi_{\alpha}(z)\,
    e^{-\frac{1}{k}(\varphi_{\alpha}(z)-\varphi_{\beta}(z))}
  \end{aligned}
\end{equation}
(as before, $\varphi_{\alpha}(z)= \alpha.\varphi(z)$ and
$\varphi_{\beta}(z)=\beta.\varphi(z)$ are the boson fields ``in the
direction'' of the corresponding screening).  With an extra boson
$\chi(z)$ added to account for the missing $\heis$, the $\hSL2_k$
algebra currents are reconstructed as
\begin{equation*}
  J^{\pm}(z)=j^{\pm}(z)\,e^{\pm\sqrt{\frac{2}{k}}\chi(z)},
  \qquad \Jnaught(z)=\sqrt{\ffrac{k}{2}}\,\partial\chi(z).
\end{equation*}

The $\hSL2$ algebra is well known, since the ``old'' studies of the
Wakimoto bosonization, to be described as a centralizer of two
fermionic screenings ``at an angle'' to each other.\footnote{The
  Wakimoto bosonization~\cite{[W]} yields two essentially different
  three-boson realizations of $\hSL2$---the ``symmetric'' and the
  ``nonsymmetric'' ones, respectively centralizing two fermionic
  screenings and one bosonic plus one fermionic screening.  The names
  refer to the ``$\jplus\leftrightarrow\jminus$ symmetric'' structure
  of~\eqref{sl2:symm} and the ``asymmetric'' structure
  of~\eqref{sl2:nonsymm}.  Somewhat broader, the ``variously
  symmetric'' realizations are discussed in~\cite{[W2n]}.}  In this
item in the list, we see again that imposing
relations~\eqref{cartan-test-log} \textit{implies} that both $-1$ in
the braiding matrix translate exactly into true fermionic screenings.

For $p\geq 3$, the Nichols algebra is given by~\cite{[Hel]}
\begin{equation*}
  \BX=T(X)/ \bigl( F_1^{2},\ [F_1,F_2]^{p},\ F_2^{2}\bigr) 
\end{equation*}
with $\dim\BX=4p$.

  \subsection*{3.2.} There are two subcases.

  \begin{description}
  \item[3.2.1] $q_{11}\in R_3$, \ $q_{12}q_{21}=q_{11}$, Cartan type
    $B_2$, $\Dynkin{\zeta}{\zeta}{-1}$\!\!,\ \ $\zeta\in R_{3}$.

    This reformulates in terms of scalar products of the screening
    momenta as
    \begin{equation*}
      \alpha.\alpha=\ffrac{2s}{3} + 2\ell\ \ (\ell\in\oZ),\quad
      2\alpha.\beta=\alpha.\alpha + 2n\ \ (n\in\oZ),
    \end{equation*}
    where $s$ is coprime with~$3$.  The braiding matrix, which is then
    parameterized as {\footnotesize$\begin{pmatrix} e^{\frac{2 i \pi s}{3}} &
      (-1)^{n+\ell} e^{\frac{i \pi  s}{3}} \\
      (-1)^{n+\ell} e^{\frac{i \pi s}{3}} & -1
    \end{pmatrix}$}, has the Cartan type $B_2$.  Its nontrivial Weyl
    reflection is {\footnotesize$\begin{pmatrix} -e^{\frac{4 i \pi s}{3}} &
      -(-1)^{n+\ell} e^{-\frac{1}{3} i \pi s}
      \\
      -(-1)^{n+\ell} e^{-\frac{1}{3} i \pi s} & -1
    \end{pmatrix}$}.

    Once again, conditions~\eqref{cartan-test-log} ensure that the
    tentative fermionic screening is such indeed, i.e., $m=0$: the
    conditions can be solved only if $(m = 0, s = 1 - 3 \ell)$ or $(m
    = 0, s = -n - 3 \ell)$.  These cases respectively yield the
    unilluminating central charges $2-\frac{6 (12 n-7)}{9 n^2+6 n-5}$
    and $-1-\frac{36}{2 n+3}+\frac{18}{n}$.

    \medskip

  \item[3.2.2\ \ (\Hel{5.11}{4})] $q_{11}\in R_3$, \
    $q_{12}q_{21}=-q_{11}$, Cartan type $B_2$,
    $\Dynkin{\zeta}{-\zeta}{-1}$\!\!,\ \ $\zeta\in R_{3}$.

    Then
    \begin{equation*}
      \alpha.\alpha=\ffrac{2s}{3},\quad
      2\alpha.\beta-\alpha.\alpha=1+2n,
    \end{equation*}
    where $s$ is coprime with~$3$.  The braiding matrix
    {\footnotesize$\begin{pmatrix} e^{\frac{2 i \pi s}{3}} &
      (-1)^{n+\ell} i e^{\frac{i\pi s}{3}}\\
      (-1)^{n+\ell} i e^{\frac{i\pi s}{3}} & -1
    \end{pmatrix}$} is of Cartan type $B_2$.  Its Weyl reflections are
  {\footnotesize$\begin{pmatrix} e^{\frac{2 i \pi s}{3}} &
    -(-1)^{n+\ell} i e^{\frac{i\pi s}{3}}\\
    -(-1)^{n+\ell} i e^{\frac{i\pi s}{3}} & -1
  \end{pmatrix}$} and\\ {\footnotesize$\begin{pmatrix} e^{\frac{4 i \pi s}{3}}
    & -(-1)^{n+\ell} e^{-i \pi (\frac{s}{3}+\frac{1}{2})}
    \\
    -(-1)^{n+\ell} e^{-i \pi (\frac{s}{3}+\frac{1}{2})} &
    -1 \end{pmatrix}$}.

    Conditions~\eqref{cartan-test-log} are solved only if $(m = 0, s =
    1 - 3 \ell)$, which leaves us with another incomprehensible
    $c=2-\frac{24 (12 n-1)}{36 n^2+60 n+1}$ (which is $c=26$ at $n=0$,
    however).

    The Nichols algebra is given by the quotient~\cite{[Hel]}
    \begin{align*}
      \BX=T(X)/\bigl( [F_1,F_1,F_2,F_1,F_2],\ F_1^{3},\ F_2^{2} \bigr)
    \end{align*}
    with $\dim\BX=36$.
  \end{description}

  \medskip

  The remaining subcases are equally unsuggestive, and the details are
  omitted.

  \subsection*{3.3} $q_0:=q_{11}q_{12}q_{21}\in R_{12}$, \
  $q_{11}=q_0^4$.  This translates into
  \begin{equation*}
    \alpha.\alpha+2\alpha.\beta=\ffrac{2r}{12} + 2j,\quad
    \alpha.\alpha=4\alpha.\alpha+8\alpha.\beta + 2n.
  \end{equation*}

  \subsection*{3.4} $q_{12}q_{21}\in R_{12}$, \
  $q_{11}=-(q_{12}q_{21})^2$, or
  \begin{equation*}
    2\alpha.\beta=\ffrac{2r}{12} + 2j,\quad
    \alpha.\alpha=4\alpha.\beta+1+2n.
  \end{equation*}

  \subsection*{3.5} $q_{12}q_{21}\in R_9$, \
  $q_{11}=(q_{12}q_{21})^{-3}$, or
  \begin{equation*}
    2\alpha.\beta=\ffrac{2r}{9} + 2j,\quad
    \alpha\alpha=-6\alpha.\beta+2n.
  \end{equation*}

  \subsection*{3.6} $q_{12}q_{21}\in R_{24}$, \
  $q_{11}=-(q_{12}q_{21})^4$, or
  \begin{equation*}
    2\alpha.\beta=\ffrac{2r}{24} + 2j,\quad
    \alpha.\alpha=8\alpha.\beta+2n.
  \end{equation*}

  \subsection*{3.7} $q_{12}q_{21}\in R_{30}$, \
  $q_{11}=-(q_{12}q_{21})^5$, or
  \begin{equation*}
    2\alpha.\beta=\ffrac{2r}{30}+2j,\quad
    \alpha.\alpha=10\alpha.\beta+1+2n.
  \end{equation*}

\section{The list, items \textbf{4.}$\pmb{*}$}
The conditions are
\begin{equation*}
  \boxed{q_{12}q_{21}\neq1, \ \ q_{11}q_{12}q_{21}\neq1, \ \
    q_{12}q_{21}q_{22}\neq1, \ \ q_{22}=-1, \ \  q_{11}\notin R_2\cup R_3,}
\end{equation*}
plus any of the conditions in cases \textbf{4.1.}--\textbf{4.8.}

In terms of the screening momenta, the common condition is
\begin{equation*}
  \beta.\beta=1 + 2m,
\end{equation*}
showing that $F_{\beta}$ is a candidate for a fermionic
screening.

\subsection*{4.1 \ (\Hel{5.11}{2})} $q_{11}\in \bigcup_{a=5}^\infty
R_a$, \ $q_{12}q_{21}=q_{11}^{-2}$, \ Cartan type $B_2$,
$\Dynkin{q}{q^{-2}}{-1}$\!\!.

In terms of screenings,
\begin{equation*}
  \alpha.\alpha=\ffrac{2}{p} + 2j,\ \ |p|\geq 5\ \ (j\in\oZ),\quad
  2\alpha.\alpha + 2\alpha.\beta = 2n.
\end{equation*}
Then the braiding matrix is parameterized as {\footnotesize$\begin{pmatrix}
    e^{\frac{2 i \pi }{p}} & (-1)^n e^{-\frac{2i \pi}{p}} \\
    (-1)^n e^{-\frac{2i \pi}{p}} & -1
  \end{pmatrix}$}.  Its Weyl reflection noncoincident with the
original matrix is
{\footnotesize$\begin{pmatrix}
    -e^{-\frac{2 i \pi }{p}} & -(-1)^n e^{\frac{2 i \pi }{p}} \\
    -(-1)^n e^{\frac{2 i \pi }{p}} & -1
  \end{pmatrix}$}.

Remarkably, once again,
conditions~\eqref{cartan-test-log} are satisfied \textit{only} for
$m=n=0$ (tentative ``peculiar'' solutions are with $|p|\leq4$).
Hence, $F_{\beta}$ is indeed a standard fermionic screening.  The
central charge is then given by
\begin{equation*}
  c = -25 + \ffrac{24}{k+3} + 6 (k+3),
\end{equation*}
where $\frac{1}{p} + j = -\frac{1}{k+1}$ (or $\frac{1}{p} + j = \half
+ \frac{1}{k+1}$).  Somewhat mysteriously, this is \textit{minus} the
central charge
\begin{equation*}
  c_{\WW{3}}=25 - \ffrac{24}{k+3} - 6 (k+3)
\end{equation*}
of the $\WW{3}$ algebra, which can be obtained by a ``partial''
Hamiltonian reduction of $\hSL3_k$ and which has a
\textit{three}-boson
realization~\cite{[Pol],[Ber],[W2n]}.\footnote{T.~Creutzig has
  suggested that this a quotient of some CFT{} that has central charge
  zero and contains the $\WW{3}$ subalgebra.  The $\WW{n}$ algebras
  can be rather versatile~\cite{[CGL],[CR]}.}

The Nichols algebra is the quotient~\cite{[Hel]}
\begin{equation*}
  \BX=T(X)\!\bigm/\!\bigl([F_1,F_1,F_1,F_2],\ F_1^{p},\
  [F_1,F_2]^{p'}, \ F_2^{2}\bigr),
\end{equation*}
where $p'=\ord(-e^{\frac{2i\pi}{p}})$, with $\dim\BX = 4 p p'$.

\bigskip

None of the remaining cases currently seems illuminating in any
respect.

\subsection*{4.2} $q_{11}\in R_5\cup R_8\cup R_{12}\cup
R_{14}\cup R_{20}$, $q_{12}q_{21}=q_{11}^{-3}$.

\subsection*{4.3} $q_{11}\in R_{10}\cup R_{18}$,
$q_{12}q_{21}=q_{11}^{-4}$.

\subsection*{4.4} $q_{11}\in R_{14}\cup R_{24}$,
$q_{12}q_{21}=q_{11}^{-5}$.

\subsection*{4.5} $q_{12}q_{21}\in R_8$, $q_{11}=(q_{12}q_{21})^{-2}$.

\subsection*{4.6} $q_{12}q_{21}\in R_{12}$,
$q_{11}=(q_{12}q_{21})^{-3}$.

\subsection*{4.7} $q_{12}q_{21}\in R_{20}$,
$q_{11}=(q_{12}q_{21})^{-4}$.

\subsection*{4.8} $q_{12}q_{21}\in R_{30}$,
$q_{11}=(q_{12}q_{21})^{-6}$.

\section{The list, items \textbf{5.}$\pmb{*}$}
I merely reproduce the original items from the list of Nichols
algebras.  The basic conditions are
\begin{equation*}
  \boxed{q_{12}q_{21}\neq1,\ \ q_{11}q_{12}q_{21}\neq1,\ \ 
    q_{12}q_{21}q_{22}\neq1,\ \  q_{11}\neq-1,\ \ q_{22}\in R_3,}
\end{equation*}
to which further conditions in \textbf{5.1.}--\textbf{5.5.} are to be
added one by one.

\subsection*{5.1} $q_0:=q_{11}q_{12}q_{21}\in R_{12}$, $q_{11}=q_0^4$,
  $q_{22}=-q_0^2$.
\subsection*{5.2} $q_{12}q_{21}\in R_{12}$, $q_{11}=q_{22}=-(q_{12}q_{21})^2$.
\subsection*{5.3} $q_{12}q_{21}\in R_{24}$, $q_{11}=(q_{12}q_{21})^{-6}$,
  $q_{22}=(q_{12}q_{21})^{-8}$.
\subsection*{5.4} $q_{11}\in R_{18}$, $q_{12}q_{21}=q_{11}^{-2}$,
  $q_{22}=-q_{11}^3$.
\subsection*{5.5} $q_{11}\in R_{30}$, $q_{12}q_{21}=q_{11}^{-3}$,
  $q_{22}=-q_{11}^5$.

\section*{$\aleph_0$.\ \ Conclusions}
Some additions to the above list (including the currently
uninteresting items?!) might hopefully follow in the future.  A
variety of isolated central charge values for lower-rank $W$-algebras
can be found in~\cite{[W-2-3]} (also see~\cite{[how]}), with
interesting possibilities of an overlap with the isolated values which
I deemed uninteresting.  I know nothing about a CFT{} counterpart of
one ``regular'' case, \textbf{2.3} (a $\mathrm{G}(3)$ reduction?).
More presentations of nonstandard type appeared recently
in~\cite{[Ag-unident]}.

As a more specific result, the construction of an octuplet algebra---a
$W_3$-counterpart of the $(1,p)$ triplet algebra---should be noted.

\medskip

It is a pleasure to acknowledge that these notes were compiled partly
as a result of the Conference on Conformal Field Theory and Tensor
Categories in Beijing.  Very special thanks go to N.~Andruskiewitsch,
J.~Fjelstad, J.~Fuchs, A.~Gainutdinov, M.~Mombelli, I.~Runkel,
C.~Schweigert, and A.~Virelizier for the very useful discussions
and, of course, to Yi-Zhi~Huang and the other organizers for their
hospitality in Beijing.  I also thank T.~Creutzig, I.~Tipunin, and
S.~Wood for useful discussions and I.~Angiono for comments, and the
referee for suggestions.  This paper was supported in part by the RFBR
grant 11-01-00830 and the RFBR--CNRS grant 09-01-93105.

\appendix
\addtocounter{section}{21}
\section{Virasoro algebra}\label{app:vir}
In CFT{}, the Virasoro algebra
\begin{equation*}
  [L_m,L_n]=(m-n) L_{m+n}+\ffrac{c}{12}(m-1)m(m+1),\quad
  m,n\in\oZ
\end{equation*}
standardly appears in the guise of an energy--momentum tensor
$T(z)=\sum_{n\in\oZ}L_n z^{-n-2}$---a (chiral) field on the complex
plane that satisfies the OPEs
\begin{equation*}
  T(z)\,T(w)=\ffrac{c/2}{(z-w)^4} + \ffrac{2T(z)}{(z-w)^2}
  + \ffrac{\partial T(z)}{z-w}.
\end{equation*}
The $c$ parameter (understood to be multiplied by the unit operator
whenever necessary) is called the central charge.

For $\theta$ bosonic fields
$\varphi(z)=(\varphi^1(z),\dots,\varphi^{\theta}(z))$ with the OPEs
\begin{equation}\label{ff-OPE}
  \varphi^i(z)\varphi^j(w)=\delta^{ij}\log(z-w),
\end{equation}
which are also frequently used in calculations in the form
\begin{equation*}
  \partial\varphi^i(z)\,\partial\varphi^j(w)=\ffrac{\delta^{ij}}{(z-w)^2},
\end{equation*}
the energy--momentum tensors are parameterized by
$\xi\in\oC^{\theta}$,
\begin{equation}\label{T-xi}
  T_{\xi}(z)=
  \half\partial\varphi(z).\partial\varphi(z)
  + \xi.\partial^2\varphi(z).
\end{equation}
The corresponding central charge is
\begin{equation}\label{c-xi}
  c_{\xi}=\theta - 12\xi.\xi \,.
\end{equation}
The OPE of $T_{\xi}(z)$ with a vertex operator $e^{\mu.\varphi(z)}$ is
\begin{equation*}
  T_{\xi}(z)\,e^{\mu.\varphi(w)}=
  \mfrac{\Delta\,e^{\mu.\varphi(w)}}{(z-w)^2}+
  \mfrac{\partial e^{\mu.\varphi(w)}}{z-w},
  \qquad \Delta = \half\, \mu.\mu -\xi.\mu\,
\end{equation*}

A screening operator is, by definition, any expression $\oint
V(\cdot)$, where $V(z)$ is a field of dimension~$\Delta=1$ (and the
contour integration is essentially equivalent to taking a residue
``after the action of $V(z)$ is evaluated'').  For $\theta=2$, any two
exponentials $e^{\alpha.\varphi(z)}$ and $e^{\beta.\varphi(z)}$ with
noncollinear $\alpha,\beta\in\oC^2$ define screening operators with
respect to the energy--momentum tensor
\begin{multline*}
  T(z)=
  \half\partial\varphi(z).\partial\varphi(z)
  \\
  - \ffrac{(2 + \alpha.\beta - \alpha.\alpha) \beta.\beta
    - 2 \alpha.\beta}{2 \delta^2}
  \,\partial^2\varphi_{\alpha}(z)  
  - \ffrac{(2 + \alpha.\beta - \beta.\beta) \alpha.\alpha
    - 2 \alpha.\beta}{2 \delta^2}
  \,\partial^2\varphi_{\beta}(z),
\end{multline*}
where $\partial\varphi_{\alpha}(z)=\alpha.\partial\varphi(z)$ and
$\partial\varphi_{\beta}(z)=\beta.\partial\varphi(z)$, and $\delta^2 =
(\alpha.\alpha)(\beta.\beta) - (\alpha.\beta)^2$.  This gives
formula~\eqref{c-rank2} for the central charge.

\medskip

Next, I show that the central charge of the $\theta$-boson
energy--momentum tensor that centralizers $\theta$ screenings $\oint
e^{\alpha_i.\varphi(z)}$, $1\leq i\leq\theta$, with linearly
independent momenta is invariant under Weyl reflections~\eqref{Rk-act}
if Eqs.~\eqref{cartan-test-log} hold.

Given the $\alpha_i$, $1\leq i\leq\theta$, the condition that all the
exponentials $e^{\alpha_i.\varphi(z)}$ have dimension~$1$ is expressed
by the system of equations for $\xi$
\begin{equation*}
  \half \alpha_i.\alpha_i - \xi.\alpha_i = 1,\qquad 1\leq i\leq\theta.
\end{equation*}
With $\xi$ written as $\xi=\sum_{j=1}^{\theta}x_j \alpha_j$, this
becomes a system for the $x_j$,
\begin{equation}\label{the-eq}
  \half \alpha_i.\alpha_i - \sum_{j=1}^{\theta}x_j \alpha_j.\alpha_i = 1,
  \quad 1\leq i\leq\theta,
\end{equation}
uniquely solvable if the $\alpha_i$ are linearly independent.

Under a Weyl groupoid operation $\mathfrak{R}^{(k)}$
in~\eqref{Rk-act}, the scalar products change and the solution $(x_j)$
also changes.    The ``old'' and ``new'' central charges are
\begin{equation*}
  c=\theta-12\sum_{\ell,j=1}^{\theta} x_{\ell} x_j\alpha_{\ell}.\alpha_j\quad
  \text{and}\quad
  \mathfrak{R}^{(k)}(c)=
  \theta-12\sum_{\ell,j=1}^{\theta}\tilde{x}_{\ell} \tilde{x}_{j}
  \mathfrak{R}^{(k)}(\alpha_{\ell}.\alpha_j), 
\end{equation*}
where the $\tilde{x}_{j}$ solve the system
``$\mathfrak{R}^{(k)}(\eqref{the-eq})$.''  With $\tilde{x}_j=x_j+y_j$,
this system becomes
\begin{multline*}
  \half\alpha_i.\alpha_i -  a_{k,i}\alpha_i.\alpha_k
  +\half a_{k,i}^2\alpha_k.\alpha_k \\
  - \sum_{j=1}^{\theta}(x_j+y_j)(\alpha_j.\alpha_i
  - a_{k,i}\alpha_j.\alpha_k- a_{k,j}\alpha_k.\alpha_i)
  + a_{k,j} a_{k,i}\alpha_k.\alpha_k)=1,\quad 1\leq i\leq\theta
\end{multline*}
(for a chosen $k$).  The claim is that this system of equations for
the ``deformation'' of the original solution is solved by the ansatz
$y_j = \delta_{j,k}y$.  To see this, substitute such $y_j$ and
use~\eqref{the-eq} in the resulting equations, which then become
\begin{equation*}
   a_{k,i}
  \Bigl(\half\alpha_k.\alpha_k( a_{k,i}+1)-\alpha_i.\alpha_k-1\Bigr)\\
  +\Bigl(y+\sum_{j=1}^{\theta}x_j a_{k,j}\Bigr)
  (\alpha_k.\alpha_i- a_{k,i}\alpha_k.\alpha_k)
  = 0,
\end{equation*}
where $1\leq i\leq\theta$.  Remarkably, if~\eqref{cartan-test-log}
holds,
then the above equations are indeed solved by
\begin{equation}\label{y-solve}
  y=1 - \ffrac{2}{\alpha_k.\alpha_k} - \sum_{j=1}^{\theta}x_j a_{k,j}.
\end{equation}

It remains to find the new central charge.  With $\tilde{x}_{j} =
x_{j} + \delta_{j,k}y$,
\begin{multline*}
  \sum_{\ell,j=1}^{\theta}\tilde{x}_{\ell} \tilde{x}_{j}
  \mathfrak{R}^{(k)}(\alpha_{\ell}.\alpha_j)={}
  \sum_{\ell,j=1}^{\theta} x_{\ell} x_j(\alpha_{\ell}.\alpha_j
  - 2 a_{k,\ell}\alpha_j.\alpha_k
  +  a_{k,l} a_{k,j}\alpha_k.\alpha_k)\\[-4pt]
  {}+ 2\!\sum_{j=1}^{\theta} y
  x_j( a_{k,j}\alpha_k.\alpha_k - \alpha_k.\alpha_j)
  + y^2\alpha_k.\alpha_k
\end{multline*}
and yet another use of~\eqref{the-eq} shows that this is
\begin{equation*}
  {}=\sum_{\ell,j=1}^{\theta} x_{\ell} x_j\alpha_{\ell}.\alpha_j
  + \Bigl(y + \sum_{j=1}^{\theta}x_j a_{k,j}\Bigr)
  \Bigl(y \alpha_k.\alpha_k + 2 - \alpha_k.\alpha_k
  + \sum_{\ell=1}^{\theta} x_{\ell} a_{k,\ell}\alpha_k.\alpha_k\Bigr),
\end{equation*}
where the last factor vanishes by virtue of~\eqref{y-solve}. \ The
central charge is invariant.

\section{$W_3$ logarithmic octuplet algebras}\label{app:W}
With the two screenings as in case \textbf{2.1} (the ``regular''
solution there, with central charge \eqref{W3-cc} and the $W(z)$
field~\eqref{W3-gen}), I propose a $W_3$ counterpart of the $(1,p)$
triplet algebra~\cite{[Kausch],[Gaberdiel-K]} by closely following the
constructions in~\cite{[FHST]}.

An \textit{octuplet} of primary fields is generated from the field
$e^{\gamma.\varphi(z)}$ with $\gamma\in\oC^2$ such that
$\gamma.\alpha^{\vee}=p$ and $\gamma.\beta^{\vee}=p$, i.e., from the
field
\begin{equation*}
  \W(z)=e^{\gamma.\varphi(z)},\qquad \gamma=\alpha^{\vee} + \beta^{\vee}
\end{equation*}
(which \textit{is} in the kernel of the two screenings $F_\alpha=\oint
e^{\alpha.\varphi}$ and $F_\beta=\oint e^{\beta.\varphi}$).  This is a
Virasoro primary field of dimension $\Delta = 3p-2$, that is,
\begin{align*}
  L_{n} \W(z)&=0,\quad n\geq 1,\\
  L_{0} \W(z)&=\Delta\W(z),\quad \Delta = 3p-2,\\
  \intertext{and, moreover, a $W_3$ primary: as is easy to verify, the
    modes of the dimension-$3$ field
    $W(z)=\sum_{n\in\oZ}W_{n}z^{-n-3}$ in~\eqref{W3-gen} act on
    $\W(z)$ such that}
  W_{n} \W(z)&=0,\quad n\geq 0.
\end{align*}
Then the long screenings~\eqref{long-scr} generate the octuplet
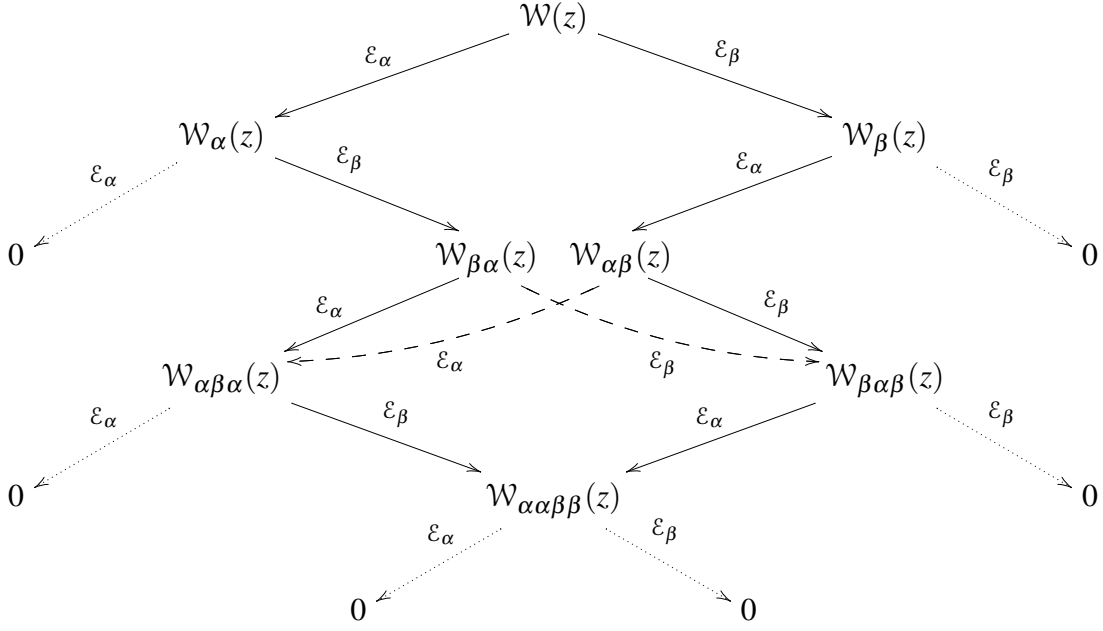
\begin{equation*}
  \xymatrix@C18pt{
    &&&&\W(z)\ar[lld]_{\Ea}\ar[rrd]^{\Eb}&&&&\\
    &&\Wa(z)\ar@{{.}{.}{>}}_{\Ea}[lld]\ar^{\Eb}[drr]+<-36pt,10pt>&&&&\Wb(z)\ar_{\Ea}[dll]+<30pt,10pt>
    \ar@{{.}{.}{>}}^{\Eb}[rrd]&&\\
    0&&&&\Wba(z)\quad \Wab(z)\ar_{\Ea}[]+<-36pt,-8pt>;[dll]
    \ar^{\Eb}[]+<36pt,-8pt>;[drr]
    \ar^{\Ea}@{{-}{--}{>}}@/^6pt/[]+<15pt,-12pt>;[dll]+<25pt,6pt>
    \ar_{\Eb}@{{-}{--}{>}}@/_6pt/[]+<-10pt,-12pt>;[drr]+<-25pt,6pt>&&&&0\\
    &&\Waba(z)\ar@{{.}{.}{>}}_{\Ea}[dll]\ar^{\Eb}[drr]&&&&\Wbab(z)\ar_{\Ea}[dll]
    \ar@{{.}{.}{>}}^{\Eb}[drr]&&\\
    0&&&&\Waabb(z)\ar@{{.}{.}{>}}_{\Ea}[dl]\ar@{{.}{.}{>}}^{\Eb}[dr]&&&&0\\
    &&&0&&0&&&
  }
\end{equation*}
Here, $\Wa(z)=\Ea\W(z)$, $\Wba(z)=\Eb\Wa(z)$, and so on, and
$\Waabb(z)=\Eb\Waba(z)=\Ea\Wbab(z)$; the dashed arrows represent maps
to the target field up to a nonzero overall factor
($\frac{(-1)^{p}}{2}$).  All the fields in the diagram are
$W_3$-algebra primaries, with the same Virasoro dimension.  All fields
below the top are of the form
$\W_{\bullet}(z)=\mathscr{P}_{\bullet}(\partial\varphi(z))\,
e^{\mu_{\bullet}.\varphi(z)}$, where the momenta $\mu_{\bullet}$ are
immediately read off from the diagram as $\mu_{\alpha}=\gamma -
\alpha^{\vee}$, \ $\mu_{\alpha\beta}=\mu_{\beta\alpha}=\gamma -
\alpha^{\vee} - \beta^{\vee} = 0$, and so on, and the
$\mathscr{P}_{\bullet}(\partial\varphi(z))$ are differential
polynomials in $\partial\varphi_{\alpha}(z)$ and
$\partial\varphi_{\beta}(z)$, of the orders
$\ord(\mathscr{P}_{\alpha})=\ord(\mathscr{P}_{\beta})=p -1 $,
$\ord(\mathscr{P}_{\alpha\beta})=\ord(\mathscr{P}_{\beta\alpha})=3 p -
2$,
$\ord(\mathscr{P}_{\alpha\beta\alpha})=\ord(\mathscr{P}_{\beta\alpha\beta})
=3 p - 3$, and $\ord(\mathscr{P}_{\alpha\alpha\beta\beta})=4 p - 4$.

Calculations in particular examples show the OPE
\begin{align*}
  \W(z)\,\Waabb(w)&=\mfrac{c_1\cdot 1}{(z-w)^{6p-4}} +
  \mfrac{c_2\,T(w)}{(z-w)^{6p-6}} +
  \mfrac{c_2/2\,\partial T(w)}{(z-w)^{6p-7}} + \dots
\end{align*}
with \textit{nonzero} coefficients (and no dimension-$3$ $W(w)$
field), and the OPEs $\Wa(z)\,\Wbab(w)$ and $\Wb(z)\,\Waba(w)$ that
start very similarly.  The adjoint-$s\ell(3)$ nature of the octuplet
manifests itself in the OPEs such as
\begin{align*}
  \Wa(z)\,\Wb(w)&=\mfrac{c_3 \W(w)}{(z-w)^{3p-2}}+\dots,\\
  \Wa(z)\,\Waba(w)&=\mathcal{O}(z-w),\\
  \Wb(z)\,\Wbab(w)&=\mathcal{O}(z-w),\\
  \Waba(z)\,\Wbab(w)&=\mfrac{c'_3 \Waabb(w)}{(z-w)^{3p-2}}+\dots.
\end{align*}


\begin{thebibliography}{99}
\bibitem{[STbr]}A.M.\;Semikhatov and I.Yu.\;Tipunin, \textit{The
    Nichols algebra of screenings}, Commun. Contemp. Math. 14 (2012)
  1250029, arXiv:1101.5810.

\bibitem{[Nich]}W.\;D.\;Nichols, \textit{Bialgebras of type one},
  Commun. Algebra 6 (1978) 1521--1552.

\bibitem{[AG]}N.\;Andruskiewitsch and M.\;Gra\~na, \textit{Braided
    Hopf algebras over non abelian finite groups}, Bol. Acad. Nacional
  de Ciencias (Cordoba) 63 (1999) 45--78 [arXiv:math$/$9802074
  [math.QA]].

\bibitem{[AS-pointed]}N.\;Andruskiewitsch and H.-J.\;Schneider,
  \textit{Pointed Hopf algebras},
  in: \textsl{New directions in Hopf algebras}, MSRI Publications 43,
  pages 1--68.  Cambridge University Press, 2002.

\bibitem{[Andr-remarks]}N.\;Andruskiewitsch, \textit{Some remarks on
    Nichols algebras}, in: \textsl{Hopf algebras}, Bergen, Catoiu and
  Chin (eds.) 25--45. M.\;Dekker (2004).

\bibitem{[AS-onthe]}N.\;Andruskiewitsch and H.-J.\;Schneider,
  \textit{On the classification of finite-dimensional pointed Hopf
    algebras}, Ann. Math. 171 (2010) 375--417 [arXiv:math$/$0502157
  [math.QA]].

\bibitem{[Heck-class]}I.\;Heckenberger, \textit{Classification of
    arithmetic root systems}, Adv. Math. 220 (2009) 59--124
  [math.QA$/$\linebreak[0]0605795].

\bibitem{[AHS]}N.\;Andruskiewitsch, I.\;Heckenberger, and
  H.-J.\;Schneider, \textit{The Nichols algebra of a semisimple
    Yetter--Drinfeld module}, Amer. J. Math., 132 (2010) 1493--1547
  [arXiv:0803.2430 [math.QA]].

\bibitem{[ARS]}N.\;Andruskiewitsch, D\;Radford, and H.-J.\;Schneider,
  \textit{Complete reducibility theorems for modules over pointed Hopf
    algebras}, J.\ Algebra, 324 (2010) 2932--2970 [arXiv:1001.3977].

\bibitem{[Heck-Weyl]}I.\;Heckenberger, \textit{The Weyl groupoid of a
    Nichols algebra of diagonal type}, Invent. Math. 164 (2006)
  175--188.

\bibitem{[GHV]}M.\;Gra\~na, I.\;Heckenberger, and L.\;Vendramin,
  \textit{Nichols algebras of group type with many quadratic
    relations}, arXiv:1004.3723.

\bibitem{[GH-lyndon]}M.\;Gra\~na and I.\;Heckenberger, \textit{On a
    factorization of graded Hopf algebras using Lyndon words}, J.
  Algebra 314 (2007) 324--343.

\bibitem{[AFGV]}N.\;Andruskiewitsch, F.\;Fantino, G.A.\;Garcia, and
  L.\;Vendramin, \textit{On Nichols algebras associated to simple
    racks}, Contemp. Math. 537, Amer. Math. Soc., Providence, RI,
  2011, pp. 31-56.  [arXiv:\linebreak[0]1006.5727].

\bibitem{[AAY]}N.\;Andruskiewitsch, I.\;Angiono, and H.\;Yamane,
  \textit{On pointed Hopf superalgebras}, Aparecer\'a en
  Contemp. Math. 544, Amer. Math. Soc., Providence, RI, 2011
  [arXiv:1009.5148].

\bibitem{[Ag-0804-standard]}I.\;Angiono, \textit{On Nichols algebras
    with standard braiding}, Algebra \& Number Theory 3 (2009) 35--106
  [arXiv:0804.0816].

\bibitem{[Ag-1008-presentation]}I.E.\;Angiono, \textit{A presentation
    by generators and relations of Nichols algebras of diagonal type
    and convex orders on root systems}, arXiv:1008.4144.

\bibitem{[Ag-1104-diagonal]}I.\;Angiono, \textit{On Nichols algebras
    of diagonal type}, arXiv:1104.0268.

\bibitem{[FHST]}J.\;Fuchs, S.\;Hwang, A.M.\;Semikhatov, and
  I.Yu.\;Tipunin, \textit{Nonsemisimple fusion algebras and the
    Verlinde formula}, Commun.\ Math.\ Phys.\ 247 (2004) 713--742
  [hep-th$/$\linebreak[0]0306274].

\bibitem{[HY]}I.\;Heckenberger and H.\;Yamane, \textit{A
    generalization of Coxeter groups, root systems, and Matsumoto's
    theorem}, Math. Z. 259 (2008) 255--276.

\bibitem{[HS]}I.\;Heckenberger and H.-J.\;Schneider, \textit{Right
    coideal subalgebras of Nichols algebras and the Duflo order on the
    Weyl groupoid}, arXiv:0909.0293 [math.QA].

\bibitem{[Heck-1+1]}I.\;Heckenberger, \textit{Finite dimensional rank
    2 Nichols algebras of diagonal type I: Examples},
  arXiv:\linebreak[0]math$/$\linebreak[0]0402350 [math.QA];
  \textit{Finite dimensional rank 2 Nichols algebras of diagonal type
    II: Classification} arXiv:\linebreak[0]math$/$\linebreak[0]0404008
  [math.QA].

\bibitem{[Hel]}M.\;Helbig, \textit{On the lifting of Nichols
    algebras}, arXiv:1003.5882 [math.QA].

\bibitem{[CH]}M.\;Cuntz and I.\;Heckenberger, \textit{Weyl groupoids
    with at most three objects}, J. Pure Appl. Algebra 213 (2009)
  1112--1128 [arXiv:0805.1810]; \textit{Finite Weyl groupoids of rank
    three}, arXiv:0912.0212 [math.QA]; \textit{Finite Weyl groupoids},
  arXiv:1008.5291.

\bibitem{[AA-gen]}N.\;Andruskiewitsch and I.E.\;Angiono, \textit{On
    Nichols algebras with generic braiding},
  arXiv:\linebreak[0]math$/$\linebreak[0]0703924, in:
  \textsl{``Modules and Comodules.''  Trends in Mathematics},
  T.\;Brzezinski, J.L.\;G\'omez Pardo, I.\;Shestakov, and P.F.\;Smith
  (eds.) (2008) 47--64.

\bibitem{[Kausch]}H.G.\;Kausch, \textit{Extended conformal algebras
    generated by a multiplet of primary fields}, Phys.\ Lett. B~259
  (1991) 448.

\bibitem{[Gaberdiel-K]}M.R.\;Gaberdiel and H.G.\;Kausch,
  \textit{Indecomposable fusion products}, Nucl.\ Phys.\ B477 (1996)
  293--318 [hep-th$/$\linebreak[0]9604026]; \textit{A rational
    logarithmic conformal field theory}, Phys.\ Lett. B~386 (1996)
  131--137 [hep-th$/$\linebreak[0]9606050]; \textit{A local
    logarithmic conformal field theory}, Nucl.\ Phys.\ B538 (1999)
  631--658 [hep-th$/$\linebreak[0]9807091].

\bibitem{[FGST]}B.L.\;Feigin, A.M.\;Gainutdinov, A.M.\;Semikhatov, and
  I.Yu.\;Tipunin, \textit{Modular group representations and fusion in
    logarithmic conformal field theories and in the quantum group
    center}, Commun. Math. Phys. 265 (2006) 47--93
  [arXiv:\linebreak[0]hep-th$/$\linebreak[0]0504093].

\bibitem{[FGST3]}B.L.\;Feigin, A.M.\;Gainutdinov, A.M.\;Semikhatov,
  and I.Yu.\;Tipunin, \textit{Logarithmic extensions of minimal
    models: characters and modular transformations}, Nucl. Phys. B757
  (2006) 303--343 [arXiv:\linebreak[0]hep-th$/$\linebreak[0]0606196].

\bibitem{[AM-3]}D.\;Adamovi\'c and A.\;Milas, \textit{On the triplet
    vertex algebra $W(p)$}, Adv.\ Math. 217 (2008) 2664--2699
  [arXiv:0707.1857v2 [math.QA]]; \textit{The $N=1$ triplet vertex
    operator superalgebras}, Commun.\ Math.\ Phys. 288 (2009) 225--270
  [arXiv:0712.0379 [math.QA]]; \textit{Lattice construction of
    logarithmic modules for certain vertex algebras},
  arXiv:\linebreak[0]0902.3417 [math.QA].

\bibitem{[FGST2]}B.L.\;Feigin, A.M.\;Gainutdinov, A.M.\;Semikhatov,
  and I.Yu.\;Tipunin, \textit{Kazhdan--Lusztig correspondence for the
    representation category of the triplet $W$-algebra in logarithmic
    CFT}, Theor. Math. Phys. 148 (2006) 1210--1235
  [arXiv:\linebreak[0]math$/$0512621 [math.QA]].

\bibitem{[FGST4]}B.L.\;Feigin, A.M.\;Gainutdinov, A.M.\;Semikhatov,
  and I.Yu.\;Tipunin, \textit{Kazhdan--Lusztig-dual quantum group for
    logarithmic extensions of Virasoro minimal models},
  J. Math. Phys. 48 (2007) 032303 [arXiv:\linebreak[0]math$/$0606506
  [math.QA]].

\bibitem{[FKW]}E.\;Frenkel, V.\;Kac, and M.\;Wakimoto,
  \textit{Characters and fusion rules for W-algebras via quantized
    Drinfeld--Sokolov reduction}, Commun. Math. Phys. 147 (1992)
  295--328.

\bibitem{[spin-4]}K.-j.\;Hamada and M.\;Takao, \textit{Spin-$4$
    current algebra}, Phys. Lett. B 209 (1988) 247--251.

\bibitem{[FfSchTh]}J.M.\;Figueroa-O'Farrill, S.\;Schrans, and
  K.\;Thielemans, \textit{On the Casimir algebra of $B_2$},
  Phys. Lett. B 263 (1991) 378--384.

\bibitem{[IT]}K.\;Ito and S.\;Terashima, \textit{Free field
    realization of $WBC_{n}$ and $WG_{2}$ algebras}, Phys. Lett. B354
  (1995) 220--231 [hep-th$/$9503165].

\bibitem{[KW]}H.G.\;Kausch and G.M.T.\;Watts, \textit{Duality in
    Quantum Toda theory and W-algebras}, Nucl. Phys. B386 (1992)
  166--192 [hep-th$/$9202070].

\bibitem{[spin-6]}J.M.\;Figueroa-O'Farrill and S.\;Schrans,
  \textit{The spin 6 extended conformal algebra}, Phys. Lett. B 245
  (1990) 471--476.

\bibitem{[Zhu]}Chuan-Jie Zhu, \textit{The complete structure of the
    $WG_2$ algebra and its BRST quantization},
  hep-th$/$\linebreak[0]9508126.

\bibitem{[BFST]}P.\;Bowcock, B.L.\;Feigin, A.M.\;Semikhatov, and
  A.\;Taormina, \textit{$\widehat{s\ell}(2|1)$ and
    $\widehat{D}(2|1;\alpha)$ as vertex operator extensions of dual
    affine $s\ell(2)$ algebras}, Commun.\ Math.\ Phys. 214 (2000)
  495--545 [hep-th$/$\linebreak[0]9907171].

\bibitem{[W2n]}B.L.\;Feigin and A.M.\;Semikhatov, \textit{$\WW{n}$
    algebras}, Nucl. Phys. B698 (2004) 409--449
  [math.QA$/$\linebreak[0]0401164].

\bibitem{[W-uni]}R.\;Blumenhagen, W.\;Eholzer, A.\;Honecker,
  K.\;Hornfeck, and R.\;Huebel, \textit{Unifying $W$-algebras},
  Phys. Lett. B332 (1994) 51--60 [hep-th$/$9404113]; \textit{Coset
    realization of unifying $W$-algebras}, Int. J. Mod. Phys. A10
  (1995) 2367--2430 [hep-th$/$9406203].

\bibitem{[W]}M.~Wakimoto, \textit{Fock representations of the affine
    Lie algebra $A_1^{(1)}$}, Commun. Math. Phys. 104 (1986) 605--609.

\bibitem{[Pol]}A.M.\;Polyakov, \textit{Gauge transformations and
    diffeomorphisms}, Int.\ J.\ Mod.\ Phys. A5 (1990) 833.

\bibitem{[Ber]}M.\;Bershadsky, \textit{Conformal field theories via
    Hamiltonian reduction}, Commun.\ Math.\ Phys. 139 (1991) 71.

\bibitem{[CGL]} T.\;Creutzig, P.\;Gao, and A.R.\;Linshaw
  \textit{Fermionic coset, critical level $W^(2)_4$-algebra and higher
    spins}, arXiv:1111.6603; \textit{A commutant realization of
    $W^(2)_n$ at critical level}, arXiv:1109.4065.

\bibitem{[CR]}T\;Creutzig and D.\;Ridout, \textit{$W$-Algebras extending
    affine $gl(1|1)$}, arXiv:1111.5049.

\bibitem{[W-2-3]}R.\;Blumenhagen, M.\;Flohr, A.\;Kliem, W.\;Nahm,
  A.\;Recknagel, and R.\;Varnhagen, \textit{$W$-Algebras with two and
    three generators}, Nucl. Phys. B 361 (1991) 255--289;
  R.\;Blumenhagen, W.\;Eholzer, A.\;Honecker, and R.\;Huebel,
  \textit{New $N=1$ extended superconformal algebras with two and
    three generators}, Int. J. Mod. Phys. A7 (1992) 7841--7871
  [hep-th$/$9207072].

\bibitem{[how]}W.\;Eholzer, A.\;Honecker, R.\;Huebel, \textit{How
    complete is the classification of $W$-symmetries?},
  Phys. Lett. B308 (1993) 42--50 [hep-th$/$9302124].

\bibitem{[Ag-unident]}I.\;Angiono, \textit{Nichols algebras of
    unidentified diagonal type}, arXiv:1108.5157.


















































\end{thebibliography}
\end{document}